\documentclass{amsart}
\usepackage{amssymb}
\usepackage{amsfonts}
\usepackage{amsmath}
\usepackage{graphicx}
\newtheorem{thm}{Theorem}[section]

\newtheorem{lemma}[thm]{Lemma}
\newtheorem{proposition}[thm]{Proposition}

\theoremstyle{definition}

\newtheorem{example}[thm]{Example}

\theoremstyle{remark}

\DeclareMathOperator{\htt}{ht} \DeclareMathOperator{\qf}{qf}
\DeclareMathOperator{\td}{t.d.}

\DeclareMathOperator{\Spec}{Spec} 
 
\numberwithin{equation}{section}

\newcommand{\field}[1]{\mathbb{#1}}
\newcommand{\C}{\field{C}}
\newcommand{\R}{\field{R}}
\newcommand{\Q}{\field{Q}}
\newcommand{\Z}{\field{Z}}

\def\1{{\rm (1)}}
\def\2{{\rm (2)}}
\def\3{{\rm (3)}}
\def\4{{\rm (4)}}
\def\5{{\rm (5)}}
\def\6{{\rm (6)}}
\def\7{{\rm (7)}}
\def\8{{\rm (8)}}
\def\9{{\rm (9)}}

\begin{document}

\title{On the Dimension Theory of Polynomial Rings over Pullbacks}

\author{S. Kabbaj}

\address{Department of Mathematics, P.O. Box 5046, King Fahd
University of Petroleum \& Minerals, Dhahran 31261, Saudi Arabia}

\email{kabbaj@kfupm.edu.sa}
\maketitle

\begin{section}{Introduction}

Since Seidenberg's (1953-54) papers \cite{S1,S2} and Jaffard's
(1960) pamphlet \cite{J} on the dimension theory of commutative
rings, the literature abounds in works exploring the prime ideal
structure of polynomial rings, including four pioneering articles by
Arnold and Gilmer on dimension sequences \cite{A,AG1,AG2,AG3}. Of
particular interest is Bastida-Gilmer's (1973) precursory article
\cite{BG} which established a formula for the Krull dimension of a
polynomial ring over a $D+M$ issued from a valuation domain.
During the last three decades, numerous papers provided in-depth
treatments of dimension theory and other related notions (such as
the S-property, strong S-property, and catenarity) in polynomial
rings over various pullback constructions. All rings considered in
this paper are assumed to be integral domains.

A polynomial ring over an arbitrary domain $R$ is subject to
Seidenberg's inequalities: $n+\dim(R)\leq\dim(R[X_{1}, ...,
X_{n}])\leq n+(n+1)\dim(R)$, $\forall\ n\geq1$. A finite-dimensional
domain $R$ is said to be Jaffard if $\dim(R[X_{1}, ..., X_{n}]) = n
+ \dim(R)$ for all $n \geq 1$; equivalently, if $\dim(R) =
\dim_{v}(R)$, where $\dim(R)$ denotes the Krull dimension of $R$ and
$\dim_{v}(R)$ its valuative dimension (i.e., the supremum of
dimensions of the valuation overrings of $R$). The study of this
class was initiated by Jaffard \cite{J}. For the convenience of the
reader, recall that, in general, for a domain $R$ with
$\dim_{v}(R)<\infty$ we have: $\dim(R)\leq \dim_{v}(R)$,
$\dim_{v}(R[X_{1}, ..., X_{n}]) = n + \dim_{v}(R)$ for all $n \geq
1$, and $\dim(R[X_{1}, ..., X_{n}]) = n + \dim_{v}(R)$ for all $n
\geq \dim_{v}(R)-1$ (Cf. \cite{ABDFK,BK,DFK,G,J}).

As the Jaffard property does not carry over to localizations (see
Example~\ref{Dim4.1} below), $R$ is said to be locally Jaffard if
$R_p$ is a Jaffard domain for each prime ideal $p$ of $R$;
equivalently, $S^{-1}R$ is a Jaffard domain for each multiplicative
subset $S$ of $R$. A locally Jaffard domain is Jaffard \cite{ABDFK}.
The class of (locally) Jaffard domains contains most classes
involved in dimension theory, including Noetherian domains
\cite{Ka}, Pr\"ufer domains \cite{G}, and universally catenarian
domains \cite{BDF1}.

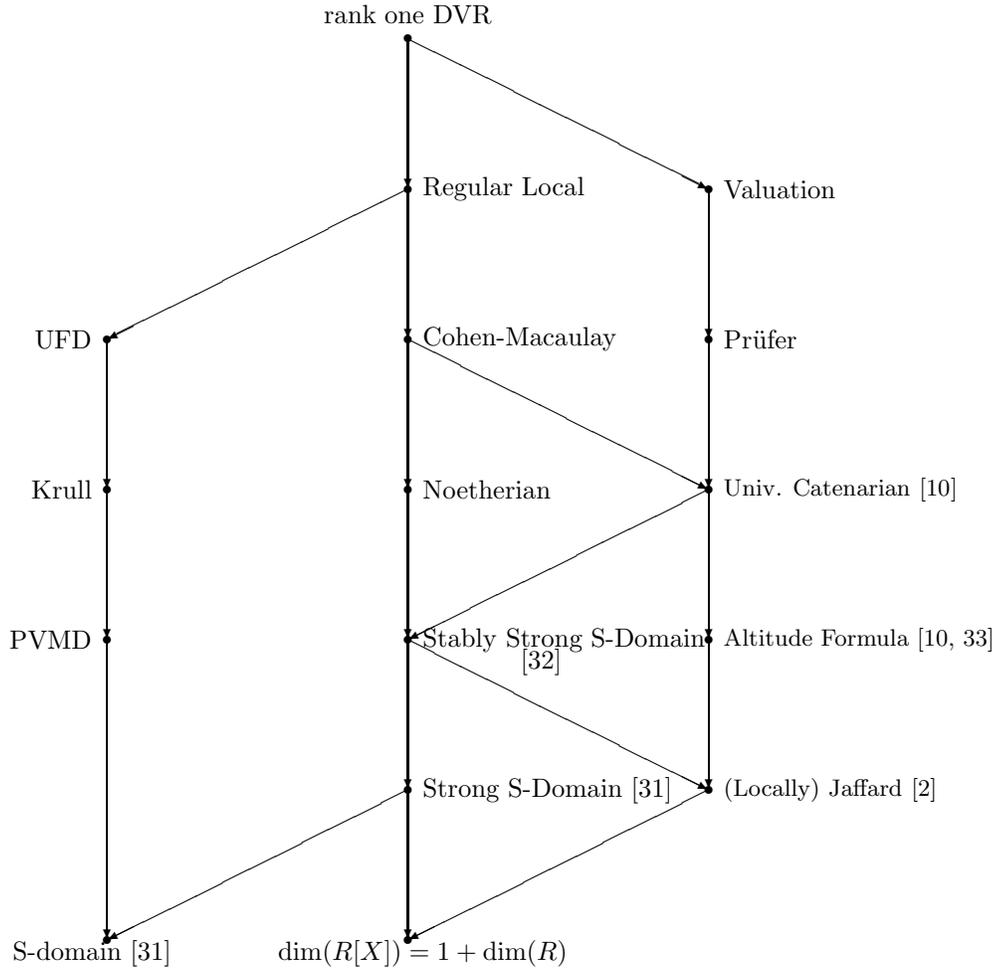
\begin{figure}[t]
\[\setlength{\unitlength}{1mm}
\begin{picture}(100,120)(0,-40)

\put(40,80){\circle*{1}}\put(40,82){\makebox(0,0)[b]{rank one DVR}}

\put(40,60){\circle*{1}}\put(42,60){\makebox(0,0)[l]{Regular Local}}
\put(80,60){\circle*{1}}\put(82,60){\makebox(0,0)[l]{Valuation}}

\put(0,40){\circle*{1}}\put(-2,40){\makebox(0,0)[r]{UFD}}
\put(40,40){\circle*{1}}\put(42,40){\makebox(0,0)[l]{Cohen-Macaulay}}
\put(80,40){\circle*{1}}\put(82,40){\makebox(0,0)[l]{Pr\"ufer}}

\put(0,20){\circle*{1}}\put(-2,20){\makebox(0,0)[r]{Krull}}
\put(40,20){\circle*{1}}\put(42,20){\makebox(0,0)[l]{Noetherian}}
\put(80,20){\circle*{1}}\put(82,20){\makebox(0,0)[l]{\small Univ.
Catenarian \cite{BDF1}}}

\put(0,0){\circle*{1}}\put(-2,0){\makebox(0,0)[r]{PVMD}}

\put(40,0){\circle*{1}}\put(42,0){\makebox(0,0)[l]{Stably Strong
S-Domain}} \put(55,-3){\makebox(0,0)[l]{\cite{MM}}}

\put(80,0){\circle*{1}}\put(82,0){\makebox(0,0)[l]{\small Altitude
Formula \cite{BDF1,M}}}

\put(40,-20){\circle*{1}}\put(42,-20){\makebox(0,0)[l]{Strong
S-Domain \cite{Ka}}}
\put(80,-20){\circle*{1}}\put(82,-20){\makebox(0,0)[l]{\small
(Locally) Jaffard \cite{ABDFK}}}

\put(0,-40){\circle*{1}}\put(-2,-40){\makebox(0,0)[t]{S-domain
\cite{Ka}}}
\put(40,-40){\circle*{1}}\put(42,-40){\makebox(0,0)[t]{$\dim(R[X])=1+\dim(R)$}}

\put(40,80){\vector(0,-2){20}} \put(40,80){\vector(2,-1){40}}

\put(40,60){\vector(0,-2){20}}\put(40,60){\vector(-2,-1){40}}
\put(80,60){\vector(0,-2){20}}

\put(0,40){\vector(0,-2){20}}
\put(40,40){\vector(0,-2){20}}\put(40,40){\vector(2,-1){40}}
\put(80,40){\vector(0,-2){20}}

\put(0,20){\vector(0,-2){20}} \put(40,20){\vector(0,-2){20}}
\put(80,20){\vector(0,-2){20}}\put(80,20){\vector(-2,-1){40}}

\put(0,0){\vector(0,-4){40}}
\put(40,0){\vector(0,-2){20}}\put(40,0){\vector(2,-1){40}}
\put(80,0){\vector(0,-2){20}}

\put(40,-20){\vector(0,-2){20}}\put(40,-20){\vector(-2,-1){40}}
\put(80,-20){\vector(-2,-1){40}}

\end{picture}\]
\caption{\tt Diagram of implications} \label{D1}
\end{figure}

In order to treat Noetherian domains and Pr\"{u}fer domains in a
unified manner, Kaplansky \cite{Ka} introduced the following
concepts: A domain $R$ is called an S-domain if, for each height-one
prime ideal $p$ of $R$, the extension $pR[X]$ in $R[X]$ has height
$1$ too; and $R$ is said to be a strong S-domain if $\frac Rp$ is an
S-domain for each prime ideal $p$ of $R$. A strong S-domain $R$
satisfies $\dim(R[X])=\dim(R)+1$. Notice that while $R[X]$ is always
an S-domain for any domain $R$ \cite{FK1}, $R[X]$ need not be a
strong S-domain even when $R$ is a strong S-domain \cite{BMRH}. Thus
$R$ is called a stably strong S-domain (also called a universally
strong S-domain) if the polynomial ring $R[X_1,...,X_n]$ is a strong
S-domain for each positive integer $n$. A stably strong S-domain is
locally Jaffard \cite{ABDFK,K1,MM}.

This review paper deals with dimension theory of polynomial rings
over certain families of pullbacks. While the literature is
plentiful, this field is still developing and many contexts are yet
to be explored. I will thus restrict the scope of the present
survey, mainly, to topics I have worked on over the last decade. The
set of pullback constructions studied includes $D+M$, $D+(X_{1},
..., X_{n})D_{S}[X_{1}, ..., X_{n}]$, $A+XB[X]$, and $D+I$.

Any unreferenced material is standard, as in \cite{B,G,J,Ka,M}. In
Figure~\ref{D1}, a diagram of implications summarizes the relations
between some spectral notions and well-known classes of integral
domains (some of which should be either finite-dimensional or
locally finite dimensional).
\end{section}

\begin{section}{Preliminaries on Pullbacks}\label{PP}
Pullbacks have proven to be useful for the construction of original
examples and counter-examples in Commutative Ring Theory. The oldest
in date is due to Krull (Cf. \cite[page 1]{BG}). However, the first
systematic investigation of a particular family of pullbacks;
namely, $D+M$ issued from valuation domains, was carried out by
Gilmer \cite[Appendix 2]{G1} and \cite{G}. Later, during the 1970s,
six ground-breaking papers \cite{BG,Gr,DP,CMZ1,BR,F} provided
further development in various pullback contexts and paved the path
for most subsequent works on these constructions. In
Figure~\ref{D2}, a diagram provides more details on the contexts
studied in these works.

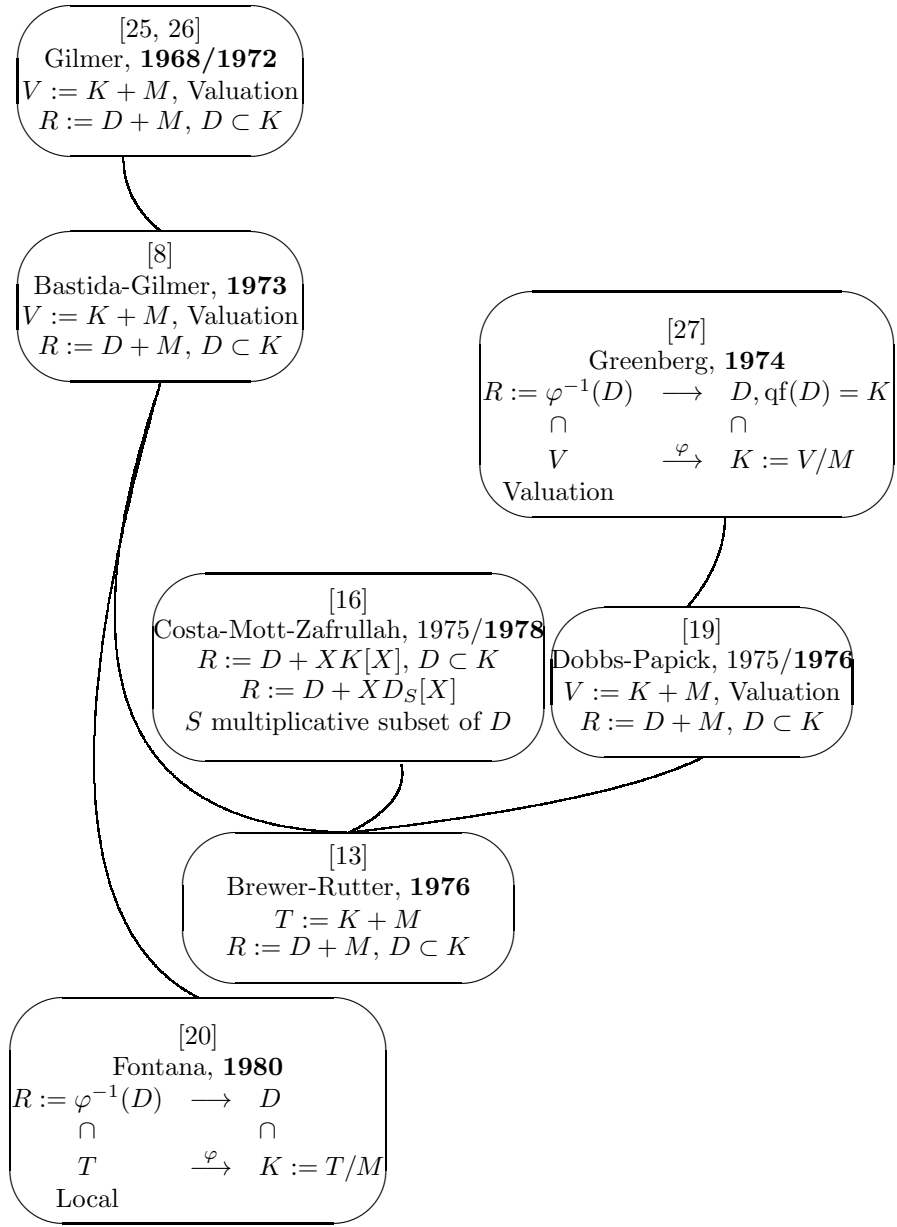
\begin{figure}[t]
\[\setlength{\unitlength}{1mm}
\begin{picture}(85,160)(0,-32)

\put(0,124){\makebox(0,0)[b]{\cite{G1,G}}}
\put(0,120){\makebox(0,0)[b]{Gilmer, {\bf 1968/1972}}}
\put(0,116){\makebox(0,0)[b]{$V:=K+M$, Valuation}}
\put(0,112){\makebox(0,0)[b]{$R:=D+M$, $D\subset K$}}
\put(0,119){\oval(38,20)}

\put(0,94){\makebox(0,0)[b]{\cite{BG}}}
\put(0,90){\makebox(0,0)[b]{Bastida-Gilmer, {\bf 1973}}}
\put(0,86){\makebox(0,0)[b]{$V:=K+M$, Valuation}}
\put(0,82){\makebox(0,0)[b]{$R:=D+M$, $D\subset K$}}
\put(0,89){\oval(38,20)}\qbezier(-5,109)(-5,103)(0,99)

\put(70,84){\makebox(0,0)[b]{\cite{Gr}}}
\put(70,80){\makebox(0,0)[b]{Greenberg, {\bf 1974}}}
\put(70,62){\makebox(0,0)[b]{$\begin{array}{ccl}
R:=\varphi^{-1}(D) & \longrightarrow                   & D, \qf(D)=K\\
\cap            &                                   & \cap\\
V               & \stackrel{\varphi}\longrightarrow    & K:=V/M\\
{\rm Valuation} &                                   &
\end{array}$}}
\put(70,76){\oval(55,30)}

\put(25,48){\makebox(0,0)[b]{\cite{CMZ1}}}
\put(25,44){\makebox(0,0)[b]{Costa-Mott-Zafrullah, 1975/{\bf 1978}}}
\put(25,40){\makebox(0,0)[b]{$R:=D+XK[X]$, $D\subset K$}}
\put(25,36){\makebox(0,0)[b]{$R:=D+XD_{S}[X]$}}
\put(25,32){\makebox(0,0)[b]{$S$ multiplicative subset of $D$}}
\put(25,41){\oval(52,25)}

\put(72,44){\makebox(0,0)[b]{\cite{DP}}}
\put(72,40){\makebox(0,0)[b]{Dobbs-Papick, 1975/{\bf 1976}}}
\put(72,36){\makebox(0,0)[b]{$V:=K+M$, Valuation}}
\put(72,32){\makebox(0,0)[b]{$R:=D+M$, $D\subset K$}}
\put(72,39){\oval(40,20)}\qbezier(75,61)(75,55)(70,49)

\put(25,14){\makebox(0,0)[b]{\cite{BR}}}
\put(25,10){\makebox(0,0)[b]{Brewer-Rutter, {\bf 1976}}}
\put(25,6){\makebox(0,0)[b]{$T:=K+M$}}
\put(25,2){\makebox(0,0)[b]{$R:=D+M$, $D\subset K$}}
\put(25,9){\oval(44,20)} \qbezier(0,79)(-20,20)(25,19)
\qbezier(32,28)(33,23)(25,19)\qbezier(72,29)(60,23)(25,19)

\put(5,-10){\makebox(0,0)[b]{\cite{F}}}
\put(5,-14){\makebox(0,0)[b]{Fontana, {\bf 1980}}}
\put(5,-32){\makebox(0,0)[b]{$\begin{array}{ccl}
R:=\varphi^{-1}(D) & \longrightarrow                   & D\\
\cap            &                                   & \cap\\
T               & \stackrel{\varphi}\longrightarrow    & K:=T/M\\
{\rm Local} &                                   &
\end{array}$}}
\put(5,-18){\oval(50,30)} \qbezier(0,79)(-20,10)(5,-3)
\end{picture}\]
\caption{\tt Diagram of various pullback contexts studied in the
1970s}\label{D2}
\end{figure}

Let's recall some results on the classical $D+M$ constructions
(i.e., those issued from valuation domains). We shall use $\qf(R)$
to denote the quotient field of a domain $R$.

\begin{thm}[\cite{G1} and \cite{DP}]\label{PP1}
Let $V$ be a valuation domain of the form $K+M$, where $K$ is a
field and $M$ is the maximal ideal of $V$. Let $D$ be a proper
subring of $K$ with $k:=\qf(D)$. Set $R:=D+M$. Then:\\
\1 $\dim(R)=\dim(V)+\dim(D)$.\\
\2 $\dim_{v}(R)=\dim(V)+\max\{\dim(W) |\ W$ is valuation on $K$
containing $D\}.$\\
\3 The integral closure of $R$ is $D'+M$, where $D'$ is the
integral closure of $D$.\\
\4 $R$ is a valuation domain  $\Leftrightarrow$  $D$ is a valuation
domain and $k=K$.\\
\5 $R$ is Pr\"ufer  $\Leftrightarrow$  $D$ is  Pr\"ufer
and $k=K$.\\
\6 $R$ is Bezout  $\Leftrightarrow$  $D$ is Bezout and
$k=K$.\\
\7 $R$ is Noetherian  $\Leftrightarrow$  $V$ is a DVR, $D=k$, and
$[K\colon k]<\infty$.\\
\8 $R$ is coherent  $\Leftrightarrow$  either ``$k=K$ and $D$ is
coherent" or ``$M$ is a finitely generated ideal of $R$." The latter
condition yields $D=k$ and $[K\colon k]<\infty$.\qed
\end{thm}

In \cite{CMZ1}, the authors established several results, similar to
the statements (1-6) and (8) above, for rings of the form $D+XK[X]$
where $K:=qf(D)$; particularly, $\dim(D+XK[X])=1+\dim(D)$ and
$\dim_{v}(D+XK[X])=1+\dim_{v}(D)$. The next result handles the
general context of $D+XD_{S}[X]$ rings.

\begin{thm}[\cite{CMZ1}]\label{PP2}
Let $D$ be an integral domain and $S$
a multiplicative subset of $D$. Set $R^{(S)}:=D+XD_{S}[X]$. Then:\\
\1 $R^{(S)}$ is GCD  $\Leftrightarrow$ $D$ is GCD and
GCD$(d,X)$ exists in $R^{(S)}, \forall\ d\in D^{*}.$\\
\2 $\dim(D_{S}[X])\leq\dim(R^{(S)})\leq\dim(D[X])$.\\
\3 If $D$ is a valuation domain, then $\dim(R^{(S)})=1+\dim(D)$.\qed
\end{thm}

in \cite{BR}, Brewer and Rutter investigated general $D+M$
constructions (i.e., issued from an integral domain not necessarily
valuation) and gave unified proofs of most results known on
classical $D+M$ and $D+XK[X]$ rings. Their result on the Krull
dimension reads as follows:

\begin{thm}[\cite{BR}]\label{PP3}
Let $T$ be an integral domain of the form $K+M$, where $K$ is a
field and $M$ is a maximal ideal of $T$. Let $D$ be a proper subring
of $K$ with $k:=\qf(D)$. Set $R:=D+M$.\\
If $k=K$, then $\dim(R)=\max\{\htt_{T}(M)+\dim(D), \dim(T)\}.$ \qed
\end{thm}

Later, Fontana \cite{F} used topological methods (particularly, his
study of amalgamated sums of
 two spectral spaces) to extend most of these results to pullbacks
(issued from local domains). We close this section by citing some
basic facts connected with the prime ideal structure of a pullback.
These will be used frequently in the sequel without explicit
mention. We shall use $\Spec(R)$ to denote the set of prime ideals
of a ring $R$.

\begin{thm}[{\cite{F} and \cite[Lemma 2.1]{ABDFK}}]\label{PP4}
Let $T$ be an integral domain, $M$ a maximal ideal of $T$, $K$ its
residue field, $\varphi:T\longrightarrow K$ the canonical
surjection, $D$ a proper subring of $K$, and $k:=\qf(D)$. Let
$R:=\varphi^{-1}(D)$ be the pullback issued from the following
diagram of canonical homomorphisms:
\[\begin{array}{ccl}
R            & \longrightarrow                 & D\\
\downarrow   &                                 & \downarrow\\
T            & \stackrel{\varphi}\longrightarrow  & K=T/M
\end{array}\]\\
\1 $M=(R:T)$ and $R/M\cong D$. \\
\2 $\Spec(R)\simeq\Spec(D)\coprod_{\Spec(K)} \Spec(T)$ (i.e.,
topological amalgamated sum)\\
\3 Assume $T$ is local. Then $M$ is a divided prime and so every
prime ideal of $R$ compares with $M$ under inclusion. If, in addition, $k=K$ then $R_{M}=T$. \\
\4 Assume $T$ is local. Then $\dim(R)=\dim(T)+\dim(D)$.\\
\5 For each prime ideal $P$ of $R$ such that $M\nsubseteq P$, there
exists a unique prime ideal $Q$ of $T$ such that $Q\cap R=P$, and
hence $T_{Q}=R_{P}$.\\
\6 For each prime ideal $P$ of $R$ such that $M\subseteq P$, there
exists a unique prime ideal $p$ of $D$ such that
$P=\varphi^{-1}(p)$, and hence $R_{P}$ can be viewed as the pullback
of $T_{M}$ and
$D_{p}$ over $K$.\\
 \7 $T$ is integral over $R$  $\Leftrightarrow$  $D=k$
and $K$ is algebraic over $k$. \qed
\end{thm}
\end{section}

\begin{section}{Dimension Theory}\label{Dim}
This section studies the Krull dimension and valuative dimension of
polynomial rings over various families of pullbacks. It also
examines the transfer of the Jaffard property to these
constructions.

In 1969, Arnold established a fundamental theorem, \cite[Theorem
5]{A}, on the dimension of a polynomial ring over an arbitrary
integral domain; namely, for any integral domain $R$ with quotient
field $K$ and for any positive integer $n$, $$\dim(R[X_{1}, ...,
X_{n}])=n+\max\{\dim(R[t_{1}, ..., t_{n}])\ |\ \{t_{i}\}_{1\leq
i\leq n}\subseteq K\}.$$ In \cite{BG}, Bastida and Gilmer generalized
this result to the case where $\{t_{i}\}_{1\leq i\leq n}$ is a
subset of an extension field of $K$. It allowed them to establish a
formula for the Krull dimension of a polynomial ring over a
classical $D+M$ as stated below:

\begin{thm}[{\cite[Theorem 5.4]{BG}}]\label{Dim1}
Let $V$ be a valuation domain of the form $K+M$, where $K$ is a
field and $M$ is the maximal ideal of $V$. Let $D$ be a proper
subring of $K$ with $k:=\qf(D)$ and let $\td(K\colon k)$ denote
the transcendence degree of $K$ over $k$. Let $n$ be a positive
integer. Set $R:=D+M$. Then:
$$\dim(R[X_{1}, ..., X_{n}])=\dim(V)+\dim(D[X_{1}, ..., X_{n}])+\min\{n,\td(K\colon k)\}.\ \ \qed$$
\end{thm}

In \cite{BK}, we refined Gilmer's statement on the valuative
dimension of a classical $D+M$ in order to build a family of
examples of Jaffard domains which are neither Noetherian nor
Pr\"ufer domains.

\begin{proposition}[{\cite[Proposition 2.1]{BK}}]\label{Dim2}
Under the same notation of Theorem~\ref{Dim1}, we have:\\
\1 $\dim_{v}(R)=\dim_{v}(D)+\dim(V)+\td(K\colon k)$.\\
\2 $R$ is a Jaffard domain  $\Leftrightarrow$  $D$ is a Jaffard
domain and $\td(K\colon k)=0$.\qed
\end{proposition}

From this result stems a first family of Jaffard domains $A_{n}$
with dimension $n+3$ which are neither Noetherian nor Pr\"ufer, for
every $n\geq1$. Indeed, the ring $B:=\Z+Y\Q(X)[Y]_{(Y)}$ is not a
Jaffard domain since $\dim(B)=2$ and $\dim_{v}(B)=3$ by
Proposition~\ref{Dim2}. For each $n\geq1$, set $A_{n}:=B[X_{1}, ...,
X_{n}]$. For $n=1$, $A_{1}=B[X_{1}]$ is a 4-dimensional Jaffard
domain, since, by Theorem~\ref{Dim1},
$\dim(B[X_{1}])=4=\dim_{v}(B)+1=\dim_{v}(B[X_{1}])$. Clearly, for
each $n\geq2$, $A_{n}$ is an $(n+3)$-dimensional Jaffard domain.
Further, $A_{1}$ is not a strong S-domain, otherwise $B$ would be so
and hence we would have $5=\dim(B[X_{1},
X_{2}])=1+\dim(B[X_{1}])=2+\dim(B)=4$, which is absurd.
Consequently, none of the rings $A_{n}$ is a strong S-domain (hence
it is neither Noetherian nor Pr\"ufer), as desired.\bigskip

We now proceed to explore a general context. Let $T$ be an integral
domain, $M$ a maximal ideal of $T$, $K$ its residue field,
$\varphi:T\longrightarrow K$ the canonical surjection, $D$ a proper
subring of $K$, and $k:=\qf(D)$. Let $R:=\varphi^{-1}(D)$ be the
pullback issued from the following diagram of canonical
homomorphisms:
\[\begin{array}{ccl}
R            & \longrightarrow                 & D\\
\downarrow   &                                 & \downarrow\\
T            & \stackrel{\varphi}\longrightarrow  & K=T/M.
\end{array}\]

\begin{thm}[{\cite[Theorem 2.6]{ABDFK}}]\label{Dim3}
Assume $T$ is local. Then:\\
\1 $\dim_{v}(R)=\dim_{v}(D)+\dim_{v}(T)+\td(K\colon k)$.\\
\2 $R$ is Jaffard  $\Leftrightarrow$ $D$ and $T$ are Jaffard and
$\td(K\colon k)=0$.\qed
\end{thm}

The next result generalizes Theorem~\ref{PP1}(1),
Theorem~\ref{PP4}(4), and Theorem~\ref{Dim3}.

\begin{thm}[{\cite[Theorem 2.11 and Corollary 2.12]{ABDFK}}]\label{Dim4}
Assume $T$ is an arbitrary domain (i.e., not necessarily local). Then:\\
\1 $\dim(R)=\max\{\dim(T), \dim(D)+\htt_{T}(M)\}$.\\
\2 $\dim_{v}(R)=\max\{\dim_{v}(T), \dim_{v}(D)+\dim_{v}(T_{M})+\td(K\colon k)\}$.\\
\3 $R$ is locally Jaffard  $\Leftrightarrow$ $D$ and $T$ are locally Jaffard and $\td(K\colon k)=0$.\\
\4 If $T$ is locally Jaffard with $\dim_{v}(T)<\infty$, $D$ is
Jaffard, and $\td(K\colon k)=0$, then $R$ is a Jaffard domain.
\qed
\end{thm}

There are examples which show that none of the hypotheses in
Theorem~\ref{Dim4}(4) is a necessary condition for $R$ to be
Jaffard. Indeed, let $V$ and $W$ be two incomparable valuation
domains of a suitable field $K$ with $n:=\dim(V)\geq 3$ and
$\dim(W)=1$. By \cite[Theorem 11.11]{N}, $T:=V\cap W$ is an
$n$-dimensional Pr\"ufer domain with two maximal ideals, say $M_{1}$
and $M$, $T_{M_{1}}=V$, and $T_{M}=W$. Let $\varphi:T\longrightarrow
T/M\cong K$ be the canonical surjection. We further require that $K$
has a subfield $k$ and a subring $D$ such that
$\dim(D)=\dim_{v}(D)=1$, $\qf(D)=k$, and $\td(K\colon k)=1$. Set
$R:=\varphi^{-1}(D)$. By Theorem~\ref{Dim4}(1) \& (2),
$\dim(R)=\dim_{v}(R)=n$. So that $R$ is Jaffard though $K$ is not
algebraic over $k$. Now, alter the above construction by taking
$n\geq 4$ and $\dim_{v}(D)=2$, so that $D$ is not Jaffard anymore,
but one can easily check that $R$ is Jaffard.

Next we proceed to the construction of the first example of a
Jaffard domain which is not locally Jaffard.

\begin{example}[{\cite[Example 3.2]{ABDFK}}]\label{Dim4.1}
Let $k$ be a field and $X_{1},X_{2},Y$ indeterminates over $k$. Set
$V_{1}:=k(X_{1},X_{2})[Y]_{(Y)}=k(X_{1},X_{2})+M_{1}$ and
$A:=k(X_{1})+M_{1}$, where $M_{1}=YV_{1}$. Let $(V,M)$ be a
one-dimensional valuation domain of the form $V=k(Y)+M$ such that
$k(Y)[X_{1},X_{2}]\subset V\subset k(X_{1},X_{2},Y)$ $\big($In order
to build such a ring, consider the valuation $v\colon
k(Y)[X_{1},X_{2}]\longrightarrow \Z^{2}$ defined by $v(X_{1})=(1,0)$
and $v(X_{2})=(0,1)$, where $\Z^{2}$ is endowed with the order
induced by the group isomorphism $i\colon \Z^{2}\longrightarrow
\Z[\sqrt{2}]$ defined by $i(a,b)=a+b\sqrt{2}\big)$. Consider the
two-dimensional valuation ring $V_{2}:=k[Y]_{(Y)}+M=k+M_{2}$ with
maximal ideal $M_{2}=Yk[Y]_{(Y)}+M$. One can easily check that
$V_{1}$ and $V_{2}$ are incomparable. By \cite[Theorem 11.11]{N},
$B:=V_{1}\cap V_{2}$ is a 2-dimensional Pr\"ufer domain with two
maximal ideals, say $N_{1}$ and $N_{2}$, $B_{N_{1}}=V_{1}$, and
$B_{N_{2}}=V_{2}$. Finally, put $R:=A\cap V_{2}$. One can show that
$R$ is semi-local with two maximal ideals ${\mathcal
M}_{1}=N_{1}\cap
 R$ and ${\mathcal M}_{2}=N_{2}\cap R$ with $R_{{\mathcal
M}_{1}}=A$ and $R_{{\mathcal M}_{2}}=V_{2}$ (Cf. \cite[Example
2.5]{DF}). Via Theorem~\ref{Dim4}, we obtain
$\dim(R)=\max\{\dim(R_{{\mathcal M}_{1}}), \dim(R_{{\mathcal
M}_{2}})\}=2$ and $\dim_{v}(R)=\max\{\dim_{v}(R_{{\mathcal M}_{1}}),
\dim_{v}(R_{{\mathcal M}_{2}})\}=2$. Thus $R$ is Jaffard but not
locally Jaffard, since $\dim(R_{{\mathcal
M}_{1}})=\dim(A)=1\not=\dim_{v}(R_{{\mathcal
M}_{1}})=\dim_{v}(A)=2$. \qed
\end{example}

The next result examines the possibility of extending
Bastida-Gilmer's result (Theorem~\ref{Dim1}) on the classical $D+M$
ring to a general context.

\begin{thm}[{\cite[Proposition 2.3 and Proposition 2.7]{ABDFK}}]\label{Dim5}
Under the same notation as above, the following statements hold.\\
\1 Assume $k=K$. Then: $\dim(R[X_{1}, ..., X_{n}])=\dim(D[X_{1},
..., X_{n}]) + \\ \dim(T[X_{1}, ..., X_{n}])-\dim(K[X_{1}, ...,
X_{n}])$, for each positive integer $n$. \\
\2 Assume $D=k$ and set $d:=\td(K\colon k)$. Then, for each
$n\geq0$, we have: $n+\dim(T)+\min\{n,d\}\leq \dim(R[X_{1}, ...,
X_{n}])\leq n+\dim_{v}(T)+d.$ \qed
\end{thm}

Now, one should design an example to show that the above
 can be strict.

\begin{example}[{\cite[Example 3.9]{ABDFK}}]\label{Dim5.1}
Let $Y_{1}, Y_{2}, U, V, Z, W$ be indeterminates over a field $k$.
Define $K:=k(Y_{1}, Y_{2})$, $S:=K(U)[V]_{(V)}$,
$R_{1}:=K(U,V,Z)[W]_{(W)}$, $A:=K(U,V)+WR_{1}$, $B:=K+VS$,
$R_{2}:=S+WR_{1}$, and $T:=K+VS+WR_{1}$. Thus, we have the following
pullbacks (with canonical homomorphisms):
\[\begin{array}{ccccc}
T           &\longrightarrow    & B             &\longrightarrow    &K\\
\downarrow  &                   & \downarrow    &                   &\downarrow\\
R_{2}       &\longrightarrow    & S             &\longrightarrow    &K(U)\\
\downarrow  &                   & \downarrow    &                   &\\
A           &\longrightarrow    & K(U,V)        &                   &\\
\downarrow  &                   & \downarrow    &                   &\\
R_{1}       &\longrightarrow    & K(U,V,Z)      &                 &
\end{array}\]
$R_{1}$ and $S$ are discrete valuation rings. Further, by applying
Theorem~\ref{PP4}(4) and Theorem~\ref{Dim3}, we obtain:

\(\begin{array}{lll}
\dim(A)=1                           &; &\dim_{v}(A)=2\\
\dim(R_{2})=\dim(S)+\dim(R_{1})=2   &; &\dim_{v}(R_{2})=3\\
\dim(B)=1                           &; &\dim_{v}(B)=2\\
\dim(T)=\dim(k)+\dim(R_{2})=2       &; &\dim_{v}(T)=4.
\end{array}\)\\
Let $\varphi:T\longrightarrow K$ be the canonical surjection and
$R:=\varphi^{-1}(k)$. The pullback $R$ has Krull dimension 2 and
valuative dimension 6. Further, $\dim(R[X])=5$ by \cite[Theorem
2.1]{F1}. Set $d:=\td(K\colon k)=2$. The desired strict
inequalities follow: $1+\dim(T)+\min\{1,d\}\lneqq \dim(R[X])\lneqq
1+\dim_{v}(T)+d.$ \qed
\end{example}

Next, we explore Costa-Mott-Zafrullah's $D+XD_{S}[X]$ construction
under a slight generalization. Let $D$ be a domain, $S$ a
multiplicative subset of $D$, and $r$ an integer $\geq 1$. Put
$R^{(S,r)}:=D+(X_{1}, ..., X_{r})D_{S}[X_{1}, ..., X_{r}]$. Let
$p\in \Spec(D)$. The $S$-coheight of $p$, denoted $S$-co$\htt(p)$,
is defined as the supremum of the lengths of all chains $p\subset
p_{1}\subset p_{2}\subset ... \subset p_{n}$ of prime ideals of $D$
with $p_{1}\cap S\not=\emptyset$. Set
$S$-$\dim(D):=\max\{S$-co$\htt(p)\ | \ p\in \Spec(D)\}$.

\begin{thm}[\cite{CMZ1} and \cite{FK1}]\label{Dim6}
Under the above notation, the following statements hold.\\
\1 $\max\{\dim(D_{S}[X_{1}, ..., X_{r}]), r+\dim(D)\}   \leq   \dim(R^{(S,r)})\\
{\hspace{2.3cm}} \leq   \min\{\dim(D[X_{1}, ..., X_{r}]),
\dim(D_{S}[X_{1}, ..., X_{r}])+S$-$\dim(D)\}.$\\
\2 $\dim_{v}(R^{(S,r)})=r+\dim_{v}(D)$.\\
\3 $D$ is Jaffard  $\Leftrightarrow$  $R^{(S,r)}$ is Jaffard and $\dim(R^{(S,r)})=r+\dim(D)$.\\
\4 $R^{(S,r)}$ is Jaffard  $\Leftrightarrow$ so is $D[X_{1}, ...,
X_{r}]$ with the same dimension as $R^{(S,r)}$. \qed
\end{thm}

Now, we provide an example to show that the Jaffard property of
$R^{(S,r)}$ does not force $D$ to be Jaffard. Here too we appeal to
pullbacks. Let $k$ be a field and $X,Y$ two indeterminates over $k$.
Put $V:= k(X)+Yk(X)[Y]_{(Y)}$ and $D:= k+Yk(X)[Y]_{(Y)}$. Clearly,
$D$ is a local domain with maximal ideal $M:=Yk(X)[Y]_{(Y)}$,
$\dim(D)=1$, and $\dim_{v}(D)=2$ by Theorem~\ref{PP1}(1) and
Proposition~\ref{Dim2}. Set $S:=D\setminus M$ and
$R^{(S,1)}:=D+XD_{S}[X]$. So $R^{(S,1)}\cong D[X]$ since $D_{M}\cong
D$. It follows that
$\dim(R^{(S,1)})=\dim(D[X])=1+\dim_{v}(D)=3=\dim_{v}(R^{(S,1)})$, as
desired.\bigskip

Next we move to a general context. Let $A\subseteq B$ an extension
of integral domains and $X$ an indeterminate over $B$. Put
$R:=A+XB[X]=\{f\in B[X]\ | \ f(0)\in A\}$. This construction was
introduced by D.D. Anderson-D.F. Anderson-Zafrullah in \cite{AAZ}.
Also, $R$ is a particular case of the constructions $B, I, D$
introduced by P.-J. Cahen \cite{C1}. Also, Int$(A)\cap B[X]=\{f\in
B[X]\ | \ f(A)\subseteq A\}$ is a subring of $R$ and hence a deeper
knowledge of $A+XB[X]$ constructions may have some interesting
impact on the integer-valued polynomial rings.

As a consequence of some general properties  of the spectrum of a
pullback \cite{F}, we state the following: First, $XB[X]$ is a prime
ideal of $R:=A+XB[X]$ with $R/XB[X]\cong A$ and hence we have an
order-isomorphism $\Spec(A)\longrightarrow\{P\in\Spec(R)\ | \
XB[X]\subseteq P\}$, $p\longmapsto p+XB[X]$. Second, $S:=\{X^{n}\ |
\ n\geq 0\}$ is a multiplicatively closed subset of $R$ and $B[X]$
with $S^{-1}R=S^{-1}B[X]=B[X, X^{-1}]$; by contraction, we obtain an
order-isomorphism $\{Q\in\Spec(B[X])\ | \ X\notin
Q\}\longrightarrow\{P\in\Spec(R)\ | \ X\notin  P\}$. Finally, the
spectral space $\Spec(R)$ is canonically homeomorphic to the
amalgamated sum of $\Spec(A)$ and $\Spec(B[X])$ over $\Spec(B)$.

For the subfamilies $D+XK[X]$ and $D+XD_{S}[X]$, it is known that
$\htt(XK[X])=\htt(XD_{S}[X])=1$. The next result probes the
situation of $XB[X]$ inside $\Spec(R)$.

\begin{thm}[{\cite[Theorem 1.2]{FIK1}}]\label{Dim7}
Let $R:=A+XB[X]$ and $N:= A\setminus \{0\}$. Then:\\
\1 $\htt_{R}(XB[X])=\dim(N^{-1}B[X])=\dim(B[X]\otimes_{A}\qf(A))$.\\
\2 $1\leq \htt_{R}(XB[X])\leq 1+\td(B\colon A)$. \qed
\end{thm}

Thus, if $\qf(A)\subseteq B$, then $\htt_{R}(XB[X])=\dim(B[X])$; and
if $A\subseteq B$ is an algebraic extension, then
$\htt_{R}(XB[X])=1$. In general, $\htt_{R}(XB[X])$ can describe all
integers between 1 and $1+\td(B\colon A)$, as shown by the
following example: Let $d$ be an integer, $t\in\{1, ..., d+1\}$, $K$
a field, and $X,X_{1}, ..., X_{d+1},Y_{1}, ..., Y_{d}$
indeterminates over $K$. Set $A:=K$ and $B:=K(X_{1}, ...,
X_{d-t+1})[Y_{1}, ..., Y_{t-1}]$. Hence $\td(B\colon A)=d$ and
$\htt_{R}(XB[X])=\dim(B[X])=t$.\bigskip

The next result studies the Krull and valuative dimensions as well
as the transfer of the Jaffard property.

\begin{thm}[{\cite[Theorems 2.1 \& 2.3]{FIK1}}]\label{Dim8}
Let $R:=A+XB[X]$ and set $k:=\qf(A)$ and $d:=\td(B\colon A)$. Then:\\
\1 $\max\{\dim(A)+\htt_{R}(XB[X]),\dim(B[X])\}\leq\dim(R)\\ {\hspace{7.3cm}}\leq\dim(A)+\dim(B[X])$.\\
\2 If $k\subseteq B$, then $\dim(R)=\dim(A)+\dim(B[X])$.\\
\3 $\dim_{v}(R)=\dim_{v}(A)+d+1$.\\
\4 $R$ is Jaffard and $\dim(R)=\dim(A)+1$ $\Leftrightarrow$ $A$
is Jaffard and $d=0$.\\
\5 If $k\subseteq B$, then: $R$ is Jaffard $\Leftrightarrow$ so is
$A$ and $\dim(B[X])=1+d$. \qed
\end{thm}

Now, one can easily construct new classes of Jaffard domains. For
instance, $\R+X\C[X,Y]$ and $\Z+X\overline{\Z}[X]$ both are
2-dimensional Jaffard domains, where $\overline{\Z}$ denotes the
integral closure of $\Z$ inside an algebraic extension of $\Q$.

The next result handles the locally Jaffard property.

\begin{thm}[{\cite[Theorems 2.8]{FIK1}}]\label{Dim9}
Let $R:=A+XB[X]$ and suppose that $A$ is a locally Jaffard domain.
Then $R$ is locally Jaffard  $\Leftrightarrow$  $B[X]$ is locally
Jaffard and $\htt_{R}(XB[X])= 1+\td(B\colon A)$. \qed
\end{thm}

We cannot knock down the hypothesis ``$A$ is locally Jaffard" to
``$A$ is Jaffard." For, assume $A$ is Jaffard but not locally
Jaffard (Example~\ref{Dim4.1}). Set $B:=\qf(A)$ and
$R:=A+XB[X]=A+X\qf(A)[X]$. In this situation $B[X]$ is locally
Jaffard and $\htt_{R}(XB[X])=1= 1+\td(B\colon A)$; whereas, $R$ is
not locally Jaffard by Theorem~\ref{Dim4}(3). Notice, however, that
the hypothesis ``$A$ is locally Jaffard" is not necessary as shown
below.

While several results concerning $D+XK[X]$ and $D+XD_{S}[X]$ are
recovered, some known results on these rings do not carry over to
the general context of $A+XB[X]$ constructions. Next, an example
provides some of these pathologies and, also, shows that the double
inequality established in Theorem~\ref{Dim8}(1) can be strict.

\begin{example}[{\cite[Example 3.1]{FIK1}}]\label{Dim10}
Let $K$ be a field and let $X, X_{1}, X_{2}, X_{3}, X_{4}$ be
indeterminates over K. Set:

\(\begin{array}{clllcll}
L   &:= &K(X_{1}, X_{2}, X_{3})     &\quad ;  &\quad V_{1}  &:=    &k+N\\
k   &:= &K(X_{1}, X_{2})            &\quad ;  &\quad D      &:=        &K(X_{1})[X_{2}]_{(X_{2})}+N\\
M   &:= &X_{4}L[X_{4}]_{(X_{4})}    &\quad ;  &\quad A      &:=        &K[X_{1}]_{(X_{1})}+M\\
N   &:= &X_{3}k[X_{3}]_{(X_{3})}    &\quad ;  &\quad B      &:=        &D+M\\
V   &:= &L+M                        &\quad ;  &\quad R  &:= &A+XB[X]
\end{array}\)

Then:\\
\1 $\max\{\dim(A)+\htt_{R}(XB[X]),\dim(B[X])\}\lneqq\dim(R)\lneqq\dim(A)+\dim(B[X])$.\\
\2 $\dim(A[X])\lneqq\dim(R)$ (in contrast with Theorem~\ref{Dim6}(1)).\\
\3 $R$ is Jaffard and $A[X]$ is not Jaffard (in contrast with Theorem~\ref{Dim6}(4)).\\
\4 $R$ is locally Jaffard and $A$ is not locally Jaffard (in
contrast with Theorem~\ref{Dim4}(3) applied to $D+XK[X]$).\medskip

Indeed, by Theorems~\ref{PP1} \&~\ref{Dim1} \&~\ref{Dim3}, $V$,
$V_{1}$, $D$, and $B$ are valuation domains of dimensions 1, 1, 2,
and 3, respectively; moreover, we have:
\begin{itemize}
\item $\dim(B[X])=\dim(B)+1=4$,
\item $\dim(A)=\dim(K[X_{1}]_{(X_{1})})+\dim(V)=2$,
\item $\dim_{v}(A)=\dim_{v}(K[X_{1}]_{(X_{1})})+\dim(V)+\td(L\colon K(X_{1}))=4$,
\item $\dim(A[X])=\dim(K[X_{1}]_{(X_{1})}[X])+\dim(V)+\min\{1,\td(L\colon K(X_{1}))\}=4$,
\item $\Spec(B)=\left\{(0) , M , P_{1}:=N+M , P_{2}:=X_{2}K(X_{1})[X_{2}]_{(X_{2})}+P_{1}\right\}$,
\item $\Spec(A)=\left\{(0) , M , Q:=X_{1}K[X_{1}]_{(X_{1})}+M\right\}$,
\item $M\cap A=P_{1}\cap A=P_{2}\cap A=M$.
\end{itemize}
Notice first that $\qf(A)=\qf(B)=\qf(V)$. Now, inside $\Spec(R)$ we
have the following chain of prime ideals (in view of the discussion
in the paragraph right before Theorem~\ref{Dim7}):
$$(0)\subsetneqq M[X]\cap R\subsetneqq P_{1}[X]\cap R\subsetneqq P_{2}[X]\cap R\subsetneqq M+XB[X]\subsetneqq Q+XB[X].$$
Therefore $\dim(R)\geq 5$, and hence $R$ is a 5-dimensional Jaffard
domain since $\dim_{v}(R)=\dim_{v}(A)+\td(B\colon A)+1=5$ by
Theorem~\ref{Dim8}. Consequently, (1) and (2) hold, and so does (3)
since $\dim_{v}(A[X])=\dim_{v}(A)+1=5$. It remains to deal with (4).
The domain $A$ is not locally Jaffard (since it is not Jaffard). Let
$P\in\Spec(R)$ with $X\notin P$. Then
$R_{P}=B[X,X^{-1}]_{PB[X,X^{-1}]}$ is a universally strong S-domain
(Cf. \cite{BDF1,MM}) and hence Jaffard (since B is a valuation
domain). So, in order to show that $R$ is locally Jaffard, it
suffices to consider the localizations with respect to the prime
ideals that contain $X$. Let $P:=p+XB[X]\in \Spec(R)$ with
$p\in\Spec(A)$. One can check that $R_{P}=A_{p}+XB[X]_{P}$ and thus
$A_{p}+XB_{p}[X]\subseteq R_{P}\subseteq A_{p}+XL[X]_{(X)}$. We
obtain, via Theorems~\ref{Dim3} \&~\ref{Dim8}, that
$\dim_{v}(R_{P})=\dim_{v}(A_{p}+XB_{p}[X])=\dim_{v}(A_{p}+XL[X]_{(X)})=\dim_{v}(A_{p})+\td(B\colon
A)+1=\dim_{v}(A_{p})+1$. We claim that $R_{P}$ is Jaffard for all
$p\in\Spec(A)$:
\begin{itemize}
\item Let $p:=(0)$. Then
$\dim(R_{P})=\htt_{R}(XB[X])=1=\dim_{v}(A_{(0)})+1$.
\item Let $p:=M$. Then the above maximal chain yields $\htt(P)=4$. Hence
$\dim(R_{P})=4=\dim_{v}(K(X_{1}))+\dim(V)+\td(L\colon
K(X_{1}))+1=\dim_{v}(A_{M})+1$. Here we view $A_{M}$ as a pullback
of $V$ and $K(X_{1})$ over $L$.
\item Let $p:=Q$. Then
$\dim(R_{P})=5=\dim_{v}(A)+1=\dim_{v}(A_{Q})+1=\dim_{v}(R_{P})$
(since $A_{Q}=A$). \qed
\end{itemize}
\end{example}

Next we move to a more general context. let $T$ be a domain, $I$ an
non-zero ideal of $T$, and $D$ a subring of $T$ such that $D\cap
I=(0)$. Throughout, $D$ will be identified with its image in $T/I$.
Also $\htt_{T}(I)$ will be assumed to be finite (though it's not
always indispensable). Let $R:=D+I$; it is a pullback determined by
the following diagram of canonical homomorphisms:
\[\begin{array}{ccl}
R:=D+I          & \longrightarrow                      & D\\
\downarrow      &                                      & \downarrow\\
T               & \longrightarrow    & T/I.
\end{array}\]
So $\Spec(R)$ is canonically homeomorphic to the amalgamated sum of
$\Spec(D)$ and $\Spec(T)$ over $\Spec(T/I)$. Precisely, $I$ is a
prime ideal of $R$ and we have the order isomorphisms:
$\Spec(D)\longrightarrow\{P\in\Spec(R)\ | \ I\subseteq P\}$,
$p\longmapsto p+I$; and $\{Q\in\Spec(T)\ | \ I\nsubseteqq
Q\}\longrightarrow\{P\in\Spec(R)\ | \ I\nsubseteqq P\}$,
$Q\longmapsto Q\cap R$.

This construction was introduced and developed by Cahen \cite{C,
C1}. Since its study has proven to be difficult in its generality,
the scope was mainly limited to the so-called $(T=B,I,D)$
almost-simple constructions (i.e., every ideal of $T$ containing $I$
is maximal). The following results -due to Cahen- approximate
$\htt_{R}(I)$ and $\dim(R)$ with respect to $\htt_{T}(I)$,
$\dim(D)$, and $\dim(T)$ in the general context.

\begin{thm}[{\cite[Proposition 5, Th\'eor\`eme 1, and Corollaire
1]{C}}]\label{Dim11}
\1 {$\htt_{T}(I)\leq \htt_{R}(I)\leq \dim(T)$.}\\
\2 $\dim(D)+\htt_{R}(I)\leq \dim(R)\leq \dim(D)+\dim(T)$.\\
\3 $\dim(R)\geq \max\{\htt_{T}(Q)+\dim(R/Q\cap R)\ |\ Q\in\Spec(T),
I\subseteq Q\}$. \qed
\end{thm}

Later, Ayache devoted his paper \cite{Ay} to the special case where
$T$ is either a finitely generated $K$-algebra or a quotient of a
power series ring in a finite number of indeterminates. He
established the following results:

\begin{thm}[\cite{Ay}]\label{Dim12}
Let $K$ be a field, $T$ a finitely generated $K$-algebra or a
quotient of a power series ring in a finite number of
indeterminates, $I$ a proper non-zero ideal of $T$, $D$ a subring of
$K$ with $k:=\qf(D)$, and $R:=D+I$. Then:\\
\1 $\dim(R)=\dim(D)+\dim(T)$.\\
\2 Assume either $T$ is a finitely generated $K$-algebra or
$\htt_{T}(I)=\dim(T)$. Then:
$\dim_{v}(R)=\dim_{v}(D)+\dim_{v}(T)+\td(K\colon k)$, and hence
$R$ is Jaffard if and only if $D$ is Jaffard and $\td(K\colon
k)=0$.  \qed
\end{thm}

We return to the general context. The next result shades more light
on $I$ within the spectrum of $R$.

\begin{lemma}[{\cite[Lemme 1.2]{FIK3}}]\label{Dim13}
Set ${\mathcal X}:=\{Q\in \Spec(T)\ |\ Q\cap R=I\}$ and ${\mathcal
Y}:=\{Q\in \Spec(T)\ |\ I\nsubseteqq Q, \exists\ Q'\in{\mathcal X} ,
(0)\subset Q\subset Q'\}$. Then:\\
\1 ${\mathcal X}\not= \emptyset$.\\
\2 ${\mathcal Y}= \emptyset$ if and only if $\htt_{R}(I)=1$.\\
\3 $\htt_{R}(I)=1+\max\{\htt_{T}(Q)\ |\ Q\in {\mathcal Y}\}$.\\
\4 If $\htt_{R[X]}(I[X])=1$, then $\td(T/Q\colon D)=0, \forall\
Q\in {\mathcal X}$. \qed
\end{lemma}

Next we show how the S-domain property is reflected on
$\htt_{R}(I)$.

\begin{thm}[{\cite[Th\'eor\`eme 1.3]{FIK3}}]\label{Dim13.1}
Assume $T$ is an S-domain. Then $R$ is an S-domain if and only if
$\htt_{R}(I)>1$ or $\td(\frac{T}{Q}\colon D)=0,\ \forall\
Q\in\Spec(T)$ such that $Q\cap R=I$. \qed
\end{thm}

In the special case where $T:=V$ is a valuation domain, one can
easily check that $\htt_{R}(I)=\htt_{V}(I)$ and
$\dim(R)=\dim(D)+\htt_{V}(I)$. Moreover, we have the following:

\begin{thm}[{\cite[Th\'eor\`eme 1.13]{FIK3}}]\label{Dim14}
Let $V$ be a valuation domain, $I$ an non-zero ideal of $V,$ $D$ a
subring of $V$ with $D\cap I=(0)$, and $R:=D+I$. Let $P_{0}$ denote
the prime ideal of $V$ that is minimal over $I$ and let $n$ be a positive integer. Then:\\
\1 $\dim_{v}(R)=\dim_{v}(D)+\dim_{v}(V_{P_{0}})+\td(\frac{V}{P_{0}}\colon D)$.\\
\2 $\dim(R[X_{1}, ..., X_{n}])=\dim(V_{P_{0}})+\dim(D[X_{1}, ..., X_{n}])+\min\{n,\td(\frac{V}{P_{0}}\colon D)\}$.\\
\3 $R$ is a Jaffard domain  $\Leftrightarrow$  $D$ is a Jaffard
domain and $\td(\frac{V}{P_{0}}\colon D)=0$. \qed
\end{thm}

Another special case is when the $D+I$ ring arises from a polynomial
ring. Namely, let $B$ be a domain, $X$ an indeterminate over $B$,
$D$ a subring of $B$, and $I$ an ideal of $B[X]$ with $I\cap B=0$.
Put $R:=D+I$. We have the following pullbacks  (with canonical
homomorphisms):

\[\begin{array}{ccl}
R:=D+I          & \longrightarrow       & D\\
\downarrow      &                       & \downarrow\\
B+I             & \longrightarrow       & B\\
\downarrow      &                       & \downarrow\\
B[X]            & \longrightarrow       & B[X]/I.
\end{array}\]

\begin{thm}[{\cite[Th\'eor\`eme 2.1]{FIK3}}]\label{Dim15}
Under the above notation, set $d:=\td(B\colon D)$. We have:\\
\1 $\dim_{v}(R)=\dim_{v}(D)+d+1$.\\
\2 $R$ is Jaffard and $\dim(R)=\dim(D)+1$ $\Leftrightarrow$ $D$ is
Jaffard and $d=0$.\qed
\end{thm}

The above result applies to the  particular context of $A+X^{n}B[X]$
constructions. Specifically, Let $A\subseteq B$ an extension of
integral domains, $X$ an indeterminate over $B$, and $n$ an integer
$\geq 1$. Put $R_{n}:=A+X^{n}B[X]$. Then
$\dim_{v}(R_{n})=\dim_{v}(A)+\td(B\colon A)+1$; and $R_{n}$ is
Jaffard and $\dim(R_{n})=\dim(A)+1$ if and only if $A$ is Jaffard
and $\td(B\colon A)=0$. Here the effect of the S-property appears
as follows: $R_{n}$ is an S-domain  if and only if
$\htt_{R_{1}}(XB[X])>1$ or $\td(B\colon A)=0$. (Since B[X] is
always an S-domain.)

In this vein, the ring $R:=\Z[(XY^{i})_{i\geq 0}]=\Z+X\Z[X,Y]$ was
shown by Ayache in \cite{Ay} to be a 3-dimensional totally Jaffard
domain \cite{C1}. In \cite{FIK3}, we improved this result by stating
that $R_{n}:=\Z[(X^{n}Y^{i})_{i\geq 0}]=\Z+X^{n}\Z[X,Y]$ is a
universally strong S-domain, for each integer $n\geq 1$.
\end{section}

\end{document}